\documentstyle[amsfonts,amscd,amssymb,latexsym,12pt]{article}

\title {\Large \bf A geometry characteristic for Banach space with $c^1$-norm }
\author{Jipu Ma }
\date{}
\begin{document}
\maketitle
\newtheorem{th}{Theorem}[section]
\newtheorem{ex}{Example}
\newtheorem{re}{Remark}[section]
\newtheorem{lem}{Lemma}[section]
\newtheorem{de}{Definition}[section]
\newtheorem{co}{Corollary}[section]
\renewcommand{\theequation}{\thesection.\arabic{equation}}
\begin{quote}
 ABSTRACT. \ Let $E$ be a Banach space and $S(E)=\{e\in E: \|e\|=1\}.$
 In this paper, a geometry characteristic for $E$ is presented by using
 a geometrical construct of $S(E).$ That is,  
  the norm of $E$ is of $c^1$  in  $ E \backslash  \{0\}$
 if and only if $S(E)$ is a  $c^1$-submanifold of $E$ with
 ${\rm codim}S(E)=1.$ The theorem is very clear, however, its proof is
 non-trivial, which shows an intrinsic connection between the continuous
 differentiability of the norm $\|\cdot\|$ in  $ E \backslash  \{0\}$  and
 differential structure of $S(E).$\\

{\bf Keywords}: Banach space, geometry, non-linear analysis, global analysis\\
{\bf 2000 Mathematics Subject  Classification}\ \ 54Exx\ \ 46Txx,\ \
58B20,\
\end{quote}
\vspace{0.3cm}
\section{ Introduction and preliminary}\label{s1}
 \vspace{0.3cm}
\indent \quad Let $E$ be a Banach space and $S(E)=\{e\in E: \|e\|=1\}.$  First of
all, we take the following example to illustrate our idea in
origin. Let $\|(x,y)\|_1=\max\{|x|,|y|\}$ and
$\|(x,y)\|=\sqrt{x^2+y^2}$ on ${\mathbb{R}}^2.$ $S({\mathbb{R}}^2)$
in the norm $\|\cdot\|_1$  is the square with length $\sqrt{2}$ of
diagonal line  and center at 0, and in the norm $\|\cdot\|$ is the unit
circle with center 0. In the second case, $S({\mathbb{R}}^2)$ is a
$c^1$- curve, but is not in the first case. This difference of
$S({\mathbb{R}}^2)$ comes from that  one of the norms is of $c^1$ in
$ {\mathbb{R}}^2\backslash \{0\}$ but the other is not. However, in
general, when $E$ is a Banach space with  $c^1$-norm in $ E
\backslash  \{0\}$, it is not known whether the geometry
structure of $S(E)$ is a characteristic for the Banach space $E$
with  $c^1$-norm in $ E\backslash \{0\}.$ In this paper, the
following theorem is proved: the norm $\|\cdot\|$ of Banach space $E$ is
of $c^1$ in $ E \backslash \{0\}$ if and only if $S(E)$ is a
$c^1$-submanifold of $E$ with ${\rm codim}S(E)=1.$ The proof is non-trivial but rather complex. Now let us recall some theorems and
definitions in global analysis,
 which are needed in the sequel.
\begin{de}\label{1.1} ([Z],[AMR],[M1-2]) Let $M$ be a topological space. $M$
 is called a $c^k$-Banach manifold ($k\geq 1$) provided that there is an
 atlas $\{(U_{\lambda} , \varphi _{\lambda },E_{\lambda})\}_{\lambda
 \in \Lambda }$ such that\\
 \quad (!`) $\bigcup _{\lambda \in \Lambda}U_{\lambda}=M.$\\
 \quad (!`!`) $\varphi _{\lambda }:U_{\lambda}\rightarrow  \varphi
 _{\lambda}
 (U_{\lambda})\subset E_{\lambda}$ is a homoemorphism where
 $E_{\lambda}$ is a Banach space.\\
\quad (!`!`!`) If $U_{\lambda}\cap U_{\mu}\neq \emptyset,$ then
$\varphi _{\lambda}\circ \varphi ^{-1}_{\mu}: E_{\mu}\rightarrow
E_{\lambda}$ and $\varphi _{\mu}\circ \varphi ^{-1}_{\lambda}:
E_{\lambda}\rightarrow E_{\mu}$ are of $c^k.$ \\The atlas
$\{(U_{\lambda} , \varphi _{\lambda },E_{\lambda})\}_{\lambda
 \in \Lambda }$ is said to be a $c^k$ differential
 structure of $M.$ \end{de}
\begin{de}\label{1.2} ([Z], [AMR], [M1-2]) Let $E$ be a Banach
space. A subset $S$ of $E$ is called a $c^k$ submanifold of $E$ if
and only if for each  $x\in S$, there exists an admissible chart
$(U, \varphi, E_{\varphi })$ of $E$ with $x\in U$ such that the following hold:\\
 \quad (!`) The chart space $E_{\varphi }$ contains a linear, closed
subspace $E_0$ which splits $E_{\varphi }.$\\
 \quad (!`!`) The chart image $\varphi (U\cap S)$ is an open set in
$E_0.$\\
 \quad (!`!`!`) $\varphi : U\rightarrow \varphi (U)\subset E_{\varphi } $
is a $c^k$-diffeomorphism.
\end{de}
(Note that if $E$ is only a $c^k$-Banach manifold $M,$ then the
condition (!`!`!`) is the same as in Definition  1.1. Here the
single chart set $\{
(E,I,E)\}$ is an atlas of $E$ and so, is simplified.)\par
 Let $E,F$ be Banach spaces, and $U$ an open set in $E.$
\begin{de}\label{1.3} ([Z], [AMR], [M4]) Suppose that
$f: x_0 \in U\subset E \rightarrow F$ is a $c^k$ map, $k\geq 1.$
$x_0$ is said to be a regular point of $f$ provided that the 
Fr$\acute{e}$chet derivative $(Df)(x_0)$ is surjective and its null space $N((Df)(x_0))$
splits $E$.
\end{de}
\begin{de}\label{1.4} ([Z], [AMR], [M4]) $y_0\in F$ is said to
be a regular value of $f$ if and only if either the preimage
$f^{-1}(y_0)$ is empty or consists of only regular points.
\end{de}
\begin{th}\label{1.1} ([Z], [AMR], [M4]) If $y_0\in F$ is a regular
value of $f,$ then the preimage $S=f^{-1}(y_0)$ is a
$c^k$-submanifold of $E$ with $T_xS=N((Df)(x)) $ for each $x\in
S.$\end{th}
\begin{th}\label{1.2}(Local normal form) ([Z]) If $f:U\subset E\rightarrow
F$ is a $c^k$-map, $k\geq 1$ and $e_0\in U,$  then there exist a
neighborhood $U_0$ at $e_0$ and $c^k$ diffeomorphism $\varphi :
U_0\rightarrow \varphi (U_0)$ with $\varphi (e_0)=0$ and $\varphi
'(e_0)=I$ such that
$$ f(e)=(Df)(e_0)\varphi (e)+ f(e_0)  \quad \forall e\in U_0.$$
\end{th}
\begin{de}\label{1.5} ([Z], [ABR]) Let $M$ be a topological space. Two charts
$(U,\varphi ,E_{\varphi })$  and $(V, \psi , E_{\psi })$are called
$c^1$-compatible if and only if $U\cap V= \emptyset,$ or $\varphi
\circ \psi ^{-1}$ from $E_{\psi}$ into $E_{\varphi }$ and  $\psi
\circ \varphi ^{-1}$ from $E_{\varphi}$ into $E_{\psi }$ are $c^1.$
\end{de}
 \indent\quad Recall that a curve $v(t)$ on the unit sphere of $E$ is
called to be of $c^1$ provided so is $(\varphi \circ v )(t)$ for an
admissible chart $(U,\varphi , E_{\varphi })$ satisfying the
conditions (!`)-(!`!`!`) in Definition  1.2;  the
 equivalent class $[v]$ generated by the curve $v$ consists of all
 $c^1$-curves $u(t)$ satisfying that $u(0)=v(0)=e$ and $(\varphi \circ v )'(0)
 =(\varphi \circ u )'(0),$ which is independent of the choice of the
 chart $(U,\varphi , E_{\varphi }).$ Let $T_eS(E)$ denote all of
 these equivalent classes, and we call it the tangent space of
 $ S(E)$ at $e\in S(E).$
\begin{re}\label{1.1}
(!`) $u(t)\in [v(t)]$ if and only if $u'(0)=v'(0)$
and $u(0)=v(0)\in S(E).$\\
(!`!`) Let $(U,\varphi ,E_{\varphi })$ be an admissible chart of $E$
at $e$ satisfying the conditions (!`)-(!`!`!`) in Definition 1.2.
Then
$$T_eS(E)=(D\varphi ^{-1})(e)E_0.$$
(Since $[v(t)]$ is determined uniquely by $v'(0)\in E,$ $T_eS(E)$ is
topology isomorphic to a closed subspace of $E,$  written as
$T_eS(E)$ still. In addition, $$(\varphi \circ v)'(0)=(D\varphi
)(e)v'(0)\in E_0,$$ and so, $T_eS(E)=(D{\varphi
}^{-1})(e)E_0.$) (For details, see [Z] and [ABR].)
\end{re}
\section{Some important lemmas and theorems  }\label{s2}
 \vspace{0.1cm}
\indent\quad In this section, the main results are as follows. A local normal
form of the $c^1$-norm $\|\cdot\|$ of Banach space $E$ is given, which
means that $\|e\|-\|e_0\|$ is locally $c^1$-diffeomorphic to a
linear functional in $E^*,$  and its proof includes a technique on
constructing the chart at each $e\in S(E)$ such that $S(E)$ becomes
a $c^1$-submanifold of $E$; if $S(E)$ is a $c^1$-submanifold of $E,$
then $P^N_{[e_0+\Delta e]}\rightarrow P^{N_0}_{[e_0]}$ as $\Delta
e\rightarrow 0,$ where $P^N_{[e_0+\Delta e]}$ and $P^{N_0}_{[e_0]}$
are the projections corresponding to the decompositions, $E=N\oplus
[e_0+\Delta e]$ and $E=N_0\oplus [e_0]$, respectively,
$N=T_{e_0+\Delta e} S(E),$ $N_0=T_{e_0} S(E),$ $[,]$ denotes the one
dimensional subspace generated by the vector in the bracket, and
$\oplus $ the topological direct sum. (This result is crucial to the
proof of the theorem. However, the result itself seems to be very
interesting and available for the study of infinite dimensional geometry.)
In order to shorten the proof of the main theorem, we  first
provide its part conclusion and preliminary theorems, lemmas as
preparations, some of which themselves are very interesting and
useful.
\begin{lem}\label{2.1} If the norm $\|\cdot\|$ of Banach space $E$ is
Fr$\acute{e}$chet differentiable at the nonzero point $e_0.$ Then
$$(D\|\cdot\|)(e_0)e_0=\|e_0\| .$$
\end{lem}
{\it Proof.}\quad  Let $\Delta e=\lambda e_0.$ Then for $\lambda $
small enough,
$$\|e+ \Delta e\|-\|e_0\|=\lambda \|e_0\|=(D\|\cdot\|)(e_0)\lambda
e_0+o(\|\Delta e\|),$$ where the term $o(\|\Delta e\|)$ means
 $\lim\limits_{\|\Delta e\|\rightarrow 0}\frac{o(\|\Delta e\|)}{\|\Delta e\|}=0.$ So
$$\|e_0\|=(D\|\cdot\|)(e_0)e_0+\lim\limits_{\lambda \rightarrow
0}\frac{o(\|\Delta e\|)}{\lambda}=(D\|\cdot\|)(e_0)e_0.$$ (note
$\lim\limits_{\lambda \rightarrow 0}\frac{o(\|\Delta e\|)}{\lambda}=
\lim\limits_{\lambda \rightarrow 0}\frac{\|e_0\|o(\|\Delta e\|)}{\|\Delta
e\|}=0$.) \hfill $\Box$\\
\begin{re}\label{2.1}
 Let $N((D\|\cdot\|)(e_0))$ denote the null space of
$(D\|\cdot\|)(e_0),$ then, under the assumption of Lemma 2.1,  $e_0\notin N((D\|\cdot\|)(e_0))$ whenever $e_0\neq0.$
\end{re}
\indent\quad The following lemma is immediate from Theorem 1.2, however, it is
 important to establish the atlas of the spheres in Banach space with $c^1$-norm.

\begin{lem}\label{2,2} Suppose that the norm $\|\cdot\|$ of Banach space
$E$ is of $c^1$ in $ E \backslash \{0\}.$
 Then for a nonzero $e_0\in E,$ there exist a
neighborhood $U_0$ at $e_0$ and a diffeomorphism $\varphi :
U_0\rightarrow \varphi (U_0)$ with $\varphi (e_0)=0$ and $\varphi
'(e_0)=I$ such that
$$ \|e\|=(D\|\cdot\|)(e_0)\varphi (e)+ \|e_0\|  \quad \forall e\in U_0.$$
\end{lem}
{\it Proof.}\quad By Remark 2.1, $(D\|\cdot\|)(e_0)e_0\neq 0$ for each
$e_0 \in E\backslash \{0\}.$ Let $N=N((D\|\cdot\|)(e_0)),$ then one has
$$ E=N((D\|\cdot\|)(e_0))\oplus [e_0] \quad \ {\rm{and}} \quad \
I=P^{[e_0]}_N+P_{[e_0]}^N .$$ Let $T^+$ denote a right inverse of
$(D\|\cdot\|)(e_0)$ such that
\begin{equation}\label{2.1}
T^+(D\|\cdot\|)(e_0)=P_{[e_0]}^N  \quad \ {\rm{and} }\quad \
(D\|\cdot\|)(e_0)T^+r=r \quad \forall r\in \mathbb{R}.
\end{equation}
(For details see [N].)  Let
\begin{equation}\label{2.2}
\varphi (e)=P^{[e_0]}_N e+T^+(\|e\|-\|e_0\|).
\end{equation}Evidently,
$\varphi (e_0)=0$
 and
$$(D\varphi )(e_0)=P^{[e_0]}_N +T^+(D\|\cdot\|)(e_0)=P^{[e_0]}_N+
P_{[e_0]}^N =I;$$ it is  clear by (2.2) that $\varphi $ is a
diffeomorphism due to the Inverse Mapping Theorem, i.e., there is a
neighborhood $U_0$ at $e_0$ such that $\varphi : U_0\rightarrow
\varphi (U_0)$ is a diffeomorphism. Note that
$(D\|\cdot\|)(e_0)P^{[e_0]}_N =0$ and $(D\|\cdot\|)(e_0)T^+(\|e\|-
\|e_0\|)=\|e\|-\|e_0\|$  by (2.1). Consequently, we obtain
$$\|e\|=(D\|\cdot\|)(e_0) \varphi (e)+\|e_0\| \quad \ \forall e\in
U_0.$$ \hfill $\Box$\par
 Suppose that the Banach space $E$ has the
decomposition $E=E_0\oplus E_1.$ Let $ e_0=P_{E_0}^{E_1}e,$
$e_1=P_{E_1}^{E_0}e,$ and $E_*=\{(e_0,e_1): \forall e_0\in E_0$
 and $e_1\in E_1\}.$  Define the norm $\|\cdot\|_*$ in $E_*$ by
 $\|(e_0,e_1)\|_*=\max\{\|e_0\|, \|e_1\|\}$ for each $(e_0,e_1)\in
 E_*.$ Evidently, $(E_*,\|\cdot\|_*)$ is a Banach space.
 The following lemma is convenient for mathematical calculus,
as one will see in the next section, although it is simple.
\begin{lem}\label {2.3} Let $\Gamma :B(E,\|\cdot\|)\rightarrow B(E_*,\|\cdot\|_*) $
be defined by $\Gamma (e)=(e_0,e_1)= (P_{E_0}^{E_1}e,
P_{E_1}^{E_0}e)$ for any $e\in E.$ Then $\Gamma \in B^{\times\
}((E,\|\cdot\|),(E_*,\|\cdot\|_*))$ and
$$\|(e_0,e_1)\|_*\leq \|\Gamma
\|\|e\| \quad \rm{and} \quad \|e\| \leq \|\Gamma
^{-1}\|\|(e_0,e_1)\|_*,
$$ where $B^{\times }((E,\|\cdot\|),(E_*,\|\cdot\|_*))$ is
the set of all invertible operators in
 $B((E,\|\cdot\|),(E_*,\|\cdot\|_*)).$ For abbreviation, write $B^{\times
 }(E,E_*)$ and $B(E,E_*)$ for them, respectively.
\end{lem}
{\it Proof.}\quad Since $N(\Gamma )=\{0\},$  $\Gamma ^{-1}(e_0,e_1)=e_0+e_1$
 for any $(e_0,e_1)\in E_*$  and $\|\Gamma \|\leq \max\{\|P^{E_0}_{E_1}\|\cdot\|P^{E_1}_{E_0}\|\},$
 the lemma is obvious  from the  Hanh-Banach Theorem.
\hfill $\Box$
\begin{th}\label {2.1} Suppose that the norm $\|\cdot\|$ of Banach space
$E$ is of $c^1$ in $E\backslash \{0\}$. Then $S=S(E)$ is a
$c^1$-submanifold of $E$ with $codim S=1.$
\end{th}
{\it Proof.}\quad By Definition  2.1, the essential to proof of
the theorem is to find an admissible chart of $E$ at each $e_0\in
S,$ fulfilling the conditions (!`)-(!`!`!`) in Definition  1.2  and
${\rm codim}E_0=1.$  By Lemma 2.2, for each $e_0\in S,$ there exist a
neighborhood $U_0$ at $e_0$ and $c^1$ diffeomorphism $\varphi :
U_0\rightarrow \varphi (U_0)\subset E ,$ $$\varphi (e)=P^{[e_0]}_N
e+T^+(\|e\|-\|e_0\|)$$ such that
$$ \|e\|=(D\|\cdot\|)(e_0)\varphi (e)+ \|e_0\|  \quad \forall e\in U_0.$$
Hereby we see  $\varphi (S\cap U_0)\subset \varphi (U_0)\cap N .$
Conversely, for any $n\in \varphi (U_0)\cap N,$ let $n=\varphi
(n_*)\in N$ for some $n_*\in U_0,$ then by the preceding equality,
$\|n_*\|=\|e_0\|=1,$ and so $n_*\in \varphi (S\cap U_0).$ This shows
$\varphi (S\cap U_0)(=\varphi (U_0)\cap N)$ is an open set in $N.$
In addition, ${\rm codim}N=\dim E/N=1$ since $(D\|\cdot\|)(e_0)e_0\neq 0.$ We
now conclude that $(U_0,\varphi , E)$ is the required chart,
 $E_0=N$ splits $E$ and ${\rm codim}E_0=\dim (E/N)=1.$
\hfill $\Box$\par
The following theorem is intuitive in geometry.
\begin{th}\label {2.2} If $S=S(E)=\{e\in E: \|e\|=1\}$ is 
a $c^1$-submanifold of Banach space $E,$ then $S_r=\{e\in E: \|e\|=r\}, r>0,$ is
also a $c^1$-submanifold of $E$ and $T_{e_0}S=T_{e_1}S_r$ for any
$e_1\in S_r, e_1=re_0,$ i.e., the tangent hyperplane $T_{e_0}S+e_0$
of $S$ at $e_0$ and $T_{e_1}S_r+e_1$ of $S_r$ at $e_1$ mutual are
parallel.
\end{th}
{\it Proof.}\quad Because of  $c^1$-submanifold $S$ of $E$ and by
Definition  1.2, for any $e_0\in S,$ there exist an admissible
chart $(U, \varphi,E_{\varphi })$ of $E$ at $e_0$ and a closed
subspace $E_0$ of $E_{\varphi }$ such that ${\rm codim}E_0=\dim
(E_{\varphi }/E_0)=1,$ $\varphi(U\cap S)$ is an open set in $E_0,$
and $\varphi :U \rightarrow \varphi(U)$ is an $c^1$-diffeomorphism.
Let $L_re=re \ \forall e\in E$ for $r\neq 0.$ Obviously, $L_r\in
B^{\times}(E).$ So both $U_1=rU$ and $r\varphi(U)$ are open sets in
$E.$ Let
$$\varphi_1(e)=(L_r\circ \varphi \circ L^{-1}_r)(e) \quad \forall
e\in E.$$  Clearly,
 $\varphi_1(U_1): U_1\rightarrow r\varphi(U)$ is a $c^1$
 diffeomorphism, and
$\varphi_1(S_r\cap U_1)=r\varphi (U\cap S)$ an open set in $E_0.$
Then by Definition  1.2, $(U_1, \varphi_1,E_{\varphi })$ is an
admissible of $E$ at any $e_1\in S_r,$ which makes that $S_r$ is a
$c^1$-submanifold of $E.$
By Remark  1.1, $T{e_0}S=(D\varphi ^{-1})(e_0)E_0$ and
$T_{e_1}S_r=(D\varphi ^{-1}_1)(e_1)E_0.$ Evidently,
$$(D\varphi_1 )(e_1)=L_r\cdot (D\varphi )(e_0)\cdot L_r^{-1},$$ from which  it follows
\begin {eqnarray*} \quad T_{e_1}S_r&=& (L_r \cdot (D\varphi^{-1})(e_0)\cdot
L^{-1}_r)E_0 =(L_r\cdot (D\varphi^{-1})(e_0))E_0 \\ &=& r
(D\varphi^{-1})(e_0)E_0=T_{e_0}S.
\end{eqnarray*}
\hfill $\Box$
\begin{lem}\label {2.4} Let $E=N_0\oplus R_0,$ then for any
closed subspace $R_1$ of $E$ satisfying $E=N_0\oplus R_1,$ there
exists an operator $\alpha \in B(R_0,N_0)$ such that
$$R_1=\{e_0+\alpha (e_0): \forall e_0\in R_0\}$$ and
$$P^{N_0}_{R_1}-P^{N_0}_{R_0}= \alpha \circ P^{N_0}_{R_0}.$$
\end{lem}
{\it Proof.}\quad Evidently,
$$P^{N_0}_{R_1}P^{N_0}_{R_0}e=
P^{N_0}_{R_1}(P^{N_0}_{R_0}e+P_{N_0}^{R_0}e)=P^{N_0}_{R_1}e=e \quad
\forall e\in R_1$$ and $$P^{N_0}_{R_0}P^{N_0}_{R_1}e=
P^{N_0}_{R_0}(P^{N_0}_{R_1}e+P_{N_0}^{R_1}e)=P^{N_0}_{R_0}e=e \quad
\forall e\in R_0,$$ i.e., ${P^{N_0}_{R_0}\mid }_{R_1}$ is the
inverse of ${P^{N_0}_{R_1}\mid }_{R_0}.$ Let $e_0=P^{N_0}_{R_0}e$
for $e\in R_1.$ Then for each $e\in R_1,$
$$e=P^{N_0}_{R_0}e+P_{N_0}^{R_0}e=e_0+P_{N_0}^{R_0}P^{N_0}_{R_1}e_0, $$
so  $\alpha ={P_{N_0}^{R_0}P^{N_0}_{R_1}\mid }_{R_0} \in B(R_0,N_0)$
such that $R_1=\{e_0+\alpha (e_0): \forall e_0\in R_0\}.$ Hereby,
$$P^{N_0}_{R_1}e=P^{N_0}_{R_0}e+\alpha (P^{N_0}_{R_0}e)=(I+\alpha
)P^{N_0}_{R_0}e \quad \forall e\in E ,$$ i.e.,
$$P^{N_0}_{R_1}- P^{N_0}_{R_0}= \alpha \circ  P^{N_0}_{R_0}.$$
\hfill $\Box$\par
The next theorem shows some interesting geometrical significance,
which is available for the study of infinite dimensional geometry.
It is also necessary for the proof of the main theorem below.
\begin{th}\label {2.3} Suppose that $S=S(E)$ is a $c^1$-submanifold
of $E$ with ${\rm codim}S=1.$ Let $N_0,N$ be the tangent spaces of
$S_{r_0}$ at $e_0$ and $S_r$ at $e_0+\Delta e,$ respectively, where
$r_0=\|e_0\|$ and $r=\|e_0+\Delta e\|.$ Then
$$P^N_{[e_0+\Delta e]}\rightarrow P^{N_0}_{[e_0]} \quad \ \rm{as} \
\Delta e\rightarrow 0.$$
\end{th}
{\it Proof.}\quad Let $(U,\varphi,E_{\varphi })$ at the point
$e_0\in E\backslash \{0\}$ be an admissible chart of $E$ satisfying
the conditions (!`)-(!`!`!`) in Definition  1.2. Assume that
$\|\Delta e\|$ is small enough such that $e_0+\Delta e\in U.$ Then
$(D\varphi ^{-1})(e) \in B^{\times }(E_\varphi ,E)$ for $e$ near
$e_0$ fulfils
$$N_0=T_{e_0}S_{r_0}=(D\varphi ^{-1})(e_0)E_0 \ \quad {\rm{and}} \quad
N =T_{e_0+\Delta e}S_r=(D\varphi ^{-1})(e_0+\Delta e)E_0.$$ Hereby
$$ N=(D\varphi ^{-1})(e_0+\Delta e)(D\varphi )(e_0)N_0.$$
 Because of ${\rm codim}S_{r_0}={\rm codim}S_r=1$ one can conclude
$$ E=N_0\oplus [e_0]=N\oplus [e_0+\Delta e].$$
Let $e_*=(D\varphi ^{-1})(e_0+\Delta e)(D\varphi )(e_0)e_0$ and
$\Psi=(D\psi ^{-1})(e_0+\Delta e)(D\varphi )(e_0).$ It then follows
 $$E=N\oplus [e_*] \quad \rm{and} \quad P^N_{[e_*]}=\Psi P^{N_0}_{[e_0]}{\Psi}^{-1}.$$
 By the continuity of $\Psi $ and ${\Psi }^{-1}$ one can assert
\begin{equation}\label{2.4 }
P^N_{[e_*]}\rightarrow P^{N_0}_{[e_0]} \quad \rm{as} \ \ \Delta
e\rightarrow 0.
\end{equation}
We now are in the position to prove the theorem. By (2.3), we have
$$P^N_{[e_*]}(e_0+\Delta e)=P^N_{[e_*]}e_0+P^N_{[e_*]} \Delta e
\rightarrow P^{N_0}_{[e_0]}e_0=e_0 \quad \rm{as} \ \  \Delta e
 \rightarrow 0.$$ By the determination of $e_*$ and the continuity of
 $\Psi $ it is immediate that $e_*\rightarrow e_0$ as $\Delta e\rightarrow 0.$
 While by Lemma  2.4, there is an operator
$\alpha \in B([e_*], N )$ such that $[e_0+\Delta e]=\{\lambda
e_*+\lambda \alpha (e_*): \forall \lambda \in \mathbb{R}\}$ so that
 $P^N_{[e_*]}(e_0+\Delta e)=\lambda e_*.$
 Then $e_0=(\lim_{\Delta e\rightarrow 0}\lambda )e_0$ and so,
 ${\lim}_{\Delta e\rightarrow 0}\lambda =1.$ We next go to show
 $ {\lim}_{\Delta e\rightarrow 0}\alpha (e_*)=0.$ Obviously, (2.3)
 implies $P_N^{[e_*]}\rightarrow P_{N_0}^{[e_0]}$ as $\Delta
 e\rightarrow 0.$ So
 $$\|P_N^{[e_*]}\Delta e\|\leq \|P_N^{[e_*]}\| \|\Delta e\|\rightarrow 0 \quad
 {\rm{and}} \quad  P_N^{[e_*]}e_0\rightarrow P_{N_0}^{[e_0]}e_0=0 $$
as $\|\Delta e\|\rightarrow 0.$  Hence it follows $\alpha
(e_*)\rightarrow 0$ as $\Delta e\rightarrow 0$ from
$P_N^{[e_*]}(e_0+\Delta e)=P_N^{[e_*]}e_0+P_N^{[e_*]}\Delta e=
\lambda \alpha (e_*)$ and $\lambda \rightarrow 1.$ In order to
complete the proof, we also need to show
$$\|P^N_{[e_0+\Delta e]}
 -P^N_{[e_*]}\|\rightarrow 0 \quad \rm{as} \ \ \Delta e\rightarrow 0 .$$ Let
$P_N^{[e_*]}h=\lambda e_*$ for any $h\in E.$ Clearly, $|\lambda
|=\frac{\|P_N^{[e_*]}h\|}{\|e_*\|}.$ Then
$$\|\alpha (P_N^{[e_*]}h)\|=|\lambda | \|\alpha (e_*)\|=\frac{ \|\alpha (e_*)\|}{\|e_*\|}
\|P_N^{[e_*]}h\|,$$ so  $$\|\alpha \|\leq \frac{ \|\alpha
(e_*)\|}{\|e_*\|}\|P_N^{[e_*]}\|.$$ In addition,  by Lemma  2.4,
 $$P^N_{[e_0+\Delta e]}-P^N_{[e_*]}=\alpha\circ P^N_{[e_*]}.$$
Thus, since $\|P^N_{[e_*]}\|\rightarrow
 \|P^{N_0}_{[e_0]}\|,$  $\|\alpha (e_*)\|\rightarrow 0,$ and
 $\|e_*\|\rightarrow \|e_0\|,$ one can assert $$\|P^N_{[e_0+\Delta e]}
 -P^N_{[e_*]}\|\rightarrow 0 \quad \rm{as} \ \ \Delta e\rightarrow 0 .$$
Finally, from $$P^N_{[e_0+\Delta
e]}-P^{N_0}_{[e_0]}=P^N_{[e_0+\Delta
e]}-P^N_{[e_*]}+P^N_{[e_*]}-P^{N_0}_{[e_0]}$$ it follows
$$P^N_{[e_0+\Delta e]}\rightarrow P^{N_0}_{[e_0]} \quad \rm{as} \ \ \Delta e\rightarrow 0.$$
\hfill $\Box$\\
\section{Main result}\label{s3}
 \vspace{0.3cm}
\begin{th}\label{3.1} Let $E$ be a Banach space. If $S=S(E)$ is a
$c^1$-submanifold of $E$ with ${\rm codim} S=1,$ then the norm $\|\cdot\|$ of
$E$ is of $c^1$ in $E\backslash \{0\}.$
\end {th}
{\it Proof.}\quad By Definition  1.2, since $S$ is a
$c^1$-submanifold of $E,$ with ${\rm codim}S=1,$ one has that for each
$e_0\in S,$ there exists a $c^1$ admissible chart $(U,\varphi,
E_{\varphi})$ of $E$ at $e_0$ such that $E_0\subset E_{\varphi}$
splits $E_{\varphi},$ $\varphi (U\cap S)$ is an open set in $E_0,$
and
$$\varphi: U\rightarrow \varphi (U) \quad \rm{and} \quad  \varphi
^{-1}: \varphi (U)\rightarrow U$$ are both  $c^1$-homoemorphisms.
 Let $\varphi (e_0)=e^0_{\varphi }\in E_0.$
  Then there exists a positive number $\eta $ such that
$$
  \quad \quad \quad \varphi ^{-1}(e^0_{\varphi }+\Delta e_{\varphi
 })-\varphi ^{-1}(e^0_{\varphi })=(D\varphi ^{-1})(e^0_{\varphi })\Delta e_{\varphi
 }+ o(\|\Delta e_{\varphi}\|)  \eqno{(3.1)} 
$$ and
$$
 e^0_{\varphi }+\Delta e_{\varphi}
 \in \varphi (U\cap S).\eqno{(3.2)} 
$$ whenever $\Delta e_{\varphi }\in
 E_0$ such that  $\|\Delta
 e_{\varphi}\|<\eta .$
Let $\tau =(D\varphi ^{-1})(e^0_{\varphi })\Delta e_{\varphi }\in
T_{e_0}S$ (see Remark 1.1). It is obvious that $\|\tau
\|<\frac{\eta}{\|(D\varphi)(e_0)\|}$ implies $\|\Delta
e_{\varphi}\|< \eta $ since $\|\Delta e_{\varphi}
\|=\|(D\varphi)(e_0) \tau\|\leq \|(D\varphi)(e_0)\| \|\tau\|$.
 Thus it follows from (3.1) and (3.2)
$$ \|e_0+\tau +o(\|\tau \|)\|-\|e_0\|=0
\quad {\rm{whenever}} \ \ \|\tau\|<\frac{\eta}{\|(D\varphi)(e_0)\|}, \eqno{(3.3)} 
$$
 (note
$\|\varphi ^{-1}(e^0_{\varphi }+\Delta e_{\varphi })\|=\|e_0\|$ by
(3.2)). Moreover,by (3.3) and the triangular inequality for the norm
$\|\cdot\|,$ it is easy to examine
$$ \|e_0+\tau\|-\|e_0\|\geq -\|o(\|\tau \|)\| \quad {\rm{and} }\quad \|o(\|\tau \|)\|\geq
\|e_0+\tau \|-\|e_0\|,$$ i.e., $\|e_0+\tau\|-\|e_0\| $ is a higher
order infinitesimal than $\|\tau\|.$ Hereby one gets
$$ \|e_0+\tau\|-\|e_0\|=o(\|\tau \|) \eqno{(3.4)} 
$$ whenever $\|\tau \|<\frac{\eta}{\|(D\varphi)(e_0)\|}.$
Hereafter, $o(\|\tau \|)$ is a real number. We claim that $e_0$ is
not in $T_{e_0}S,$ since
 otherwise it leads to the contradiction  that by (3.4)
$$ \|e_0+\lambda e_0\|-\|e_0\|=o(\|\tau \|) \quad {\rm{whenever} }\ \
|\lambda | < \frac{\eta}{\|(D\varphi)(e_0)\|}$$ but  by computing
directly,
$$\|e_0+\tau\| -\|e_0\|=\lambda \|e_0\|.$$
So, by ${\rm codim}S=1,$ one has  $E=N_0\oplus [e_0]$ where $N_0= T_{e_0}S.$
Next we show that the norm $\|.\|$ of $E$ is Fr$\acute{e}$chet
differentiable at each $e_0\in S.$ Let $h=\tau+\lambda e_0 $ for any
$h\in E$ where $\tau \in T_{e_0}S.$ By computing directly,
\begin {eqnarray*} \quad
&&\|e_0+(\tau +\lambda e_0)\|-\|e_0\|\\&=&\|e_0+\tau +\lambda e_0\|-
\|e_0+\lambda e_0\|+\|e_0+\lambda e_0\|-\|e_0\|\\&=&(1+\lambda
)\|e_0+\frac{\tau}{1+\lambda }\|-\|(1+\lambda )e_0\|+\|(1+\lambda
)e_0\|-\|e_0\|\\ &=& (1+\lambda )\{\|e_0+\frac{\tau}{1+\lambda
}\|-\|e_0\|\}+\lambda \|e_0\|
\end{eqnarray*}for $ |\lambda |<1.$
Further, applying (3.4) to $\|e_0+\frac{\tau}{1+\lambda }\|$,
$$\|e_0+(\tau +\lambda e_0)\|-\|e_0\|=\lambda \|e_0\|+  (1+\lambda )o(\|\frac{\tau}{1+\lambda }\|) $$
for $\|\frac{\tau}{1+\lambda }\|<\frac{\eta}{\|(D\varphi)(e_0)\|}$
and $|\lambda |<1.$ Since for $|\lambda |<1,$ $$\frac{\|\tau
\|}{1+\lambda}\rightarrow 0 \Leftrightarrow \|\tau \|\rightarrow 0$$
and $$\frac{(1+\lambda )o(\frac{\|\tau \|}{1+\lambda})}{{\|\tau
\|}}=\frac{o(\frac{\|\tau \|}{1+\lambda })}{\frac{\|\tau
\|}{1+\lambda }},$$ it is clear that $(1+\lambda )o(\frac{\|\tau
\|}{1+\lambda}),$ written still by $o(\|\tau \|),$ is also a higher
order infinitesimal than $\|\tau \|.$
   Therefore, one can assert
$$
\|e_0+h\|-\|e_0\|=\lambda \|e_0\|+o(\|\tau \|),\eqno{(3.5)} 
$$ whenever $\|\tau
\|<\frac{\eta}{2\|(D\varphi )(e_0)\|}$ and $|\lambda
|<\frac{1}{2}.$  Let $\delta =\min\{\frac{\eta }{2\|(D\varphi
 )(e_0)\|},\frac{1}{2}\},$ and $h=\tau +\lambda e_0$ for each $h\in E,$
 where $\tau \in T_{e_0}S.$ Define $\Gamma h= (\tau ,\lambda)$ by
 the same way as in Lemma  2.3. Then from
 $$ \|\Gamma \|\|h\|\geq \|\Gamma h\|_*=\max\{\|\tau\|,|\lambda
|\}$$ it follows that for any $h$ such that $\|h\|< \|\Gamma \|^{-1}\delta ,$
$$\|\tau \|<\frac{\eta }{2\|(D\varphi
)(e_0)\|} \ \ \rm{and} \ \ |\lambda |<\frac {1}{2}.$$
 Thus, by (3.4) we have
$$
\|e_0+h\|-\|e_0\|=\lambda \|e_0\|+o(\|\tau \|) \quad {\rm{whenever}} \
\ \|h\|<\|\Gamma \|^{-1}\delta .\eqno{(3.6)} 
$$
 In order to prove the Fr$\acute{e}$chet differentiability of the norm
$\|\cdot\|$ in $T_{e_0}S,$ we also have to show
$$\lim_{h\rightarrow 0} \frac {o(\|\tau \|)}{\|h\|}=\lim_{\tau \rightarrow 0} \frac {o(\|\tau
\|)}{\|\tau \|}.$$ By Lemma  2.3  it is easy to see
$$\|h\|\rightarrow 0\Leftrightarrow \|(\tau ,\lambda)\|_*\rightarrow
0\Rightarrow \|\tau \|\rightarrow 0,$$ and
$$\frac{o(\|\tau \|)}{\|h\|}\leq \frac{o(\|\tau \|)}{\|\Gamma\|^{-1}\|(\tau ,\lambda)\|_*}
\leq \frac{o(\|\tau \|)}{\|\Gamma\|^{-1}\|\tau \|}.$$ Then by (3.5)
$$
\|e_0+h\|-\|e_0\|=\lambda \|e_0\|+o(\|h \|) \quad{ \rm{whenever}} \ \
\|h\|<\|\Gamma \|^{-1}\delta .
$$
 Let $h=P_{N _0} ^{[e_0]}h+P^{N_0}_{[e_0]}h ,$ and
$P^N_{[e_0]}h=\lambda e_0$ for each point $e_0\in S.$ Define a
bounded linear functional $f_{e_0}\in E^*$ as follows:
$$ f_{e_0}(h)=\lambda \|e_0\| \quad \quad \forall h\in E,$$ which
satisfies
$$|f_{e_0}(h)|=\|P^N_{[e_0]}h\|\leq \|P^N_{[e_0]}\|\|h\|.$$
Finally one gets by (3.6)
$$\|e_0+h\|-\|e_0\|=f_{e_0}(h)+o(\|h\|) \quad {\rm{whenever}} \ \
\|h\|<\|\Gamma \|^{-1}\delta .$$ This proves that
$$
(D\|\cdot\|)(e_0)h=f_{e_0}(h),
$$
 i.e., the norm $\|\cdot\|$ is Fr$\acute{e}$chet
differentiable at each $e_0\in S.$\par
 Next we show that the norm
$\|\cdot\|$ is Fr$\acute{e}$chet differentiable for each $e\in
E\backslash \{0\}.$ Let $e_1=\|e_1\|e_0$ for each $e_1\in
E\backslash \{0\},$ then $\|e_0\|=1.$ Replace $e_0, \varphi , S$ and
$U$ above by $e_1,\varphi _1=r\varphi ,S_r, $ and $U_1=rU,$
respectively, where $ r=\|e_1\|.$ Note that $E_{\varphi }$ and $E_0$
keep invariant as shown in the proof of Theorem  2.2. Repeat the
process above.
Then, the following results follow in turn,\\
(!`) there is a positive number $\eta $ such that
$$ \varphi_1 ^{-1}(e^1_{\varphi_1 }+\Delta e_{\varphi
 })-\varphi_1 ^{-1}(e^1_{\varphi_1 })=(D\varphi_1 ^{-1})(e^1_{\varphi_1 })\Delta
 e_{\varphi
 }+ o(\|\Delta e_{\varphi}\|)$$
 whenever  $\|\Delta e_{\varphi}\|<\eta $ and $ \Delta e_{\varphi}\in E_0.$\\
(!`!`) let $h=\tau+\lambda e_1$ for any $h\in E,$ where $\tau \in
T_{e_1}S_r,$  $\delta =min \{\frac{\eta }{2(D\varphi_1 }
 , \frac{\|e_1\|}{2}\},$ and $\Gamma (h)=(\tau ,\lambda e_1),$ then  $$\|e_1+h\|-\|e_1\|=\lambda \|e_1\|+o(\|h \|) \quad \rm{whenever} \ \
\|h\|<\|\Gamma \|^{-1}\delta .$$\\
(!`!`!`) for  any $h\in E,$ let $h=P_N^{[e_1]}h+P^N_{[e_1]}h ,$ and
$P^N_{[e_1]}h=\lambda e_1$ (where $N=T_{e_0}S=T_{e_1}S_r$ by Theorem
2.2), then
$$(D\|\cdot\|)(e_1)h=f_{e_1}(h),$$
where $f_{e_1}$ is the bounded linear functional determined by
$f_{e_1}(h)=\lambda \|e_1\| $ as $P^N_{[e_1]}h=\lambda e_1.$\par
 To the end of the proof, it remains to examine the continuity of $(D\|\cdot\|).$
Let $(U,\varphi ,E_{\varphi })$ at any point $e_0\in E\backslash
\{0\}$ be an admissible chart of $E$ satisfying the conditions
(!`)-(!`!`!`) in Definition  1.2. Assume that $\|\Delta e\|$ is
small enough such that $e_0+\Delta e\in U.$ Let $r_0=\|e_0\|, \
r=\|e_0+\Delta e\|.$ Thus
$$N_0=T_{e_0}S_{r_0}=(D\varphi ^{-1})(e_0)E_0 \ \quad \rm{and} \quad
N =T_{e_0+\Delta e}S_r=(D\varphi ^{-1})(e_0+\Delta e)E_0.$$ Hereby
$$ N=(D\varphi ^{-1})(e_0+\Delta e)(D\varphi )(e_0)N_0.$$
 Because of ${\rm codim}S_{r_0}={\rm codim}S_r=1$ one can conclude
$$ E=N_0\oplus [e_0]=N\oplus [e_0+\Delta e].$$
By Theorem  2.3,
$$P^N_{[e_0+\Delta e]}\rightarrow P^{N_0}_{[e_0]} \quad \ \rm{as} \
\Delta e\rightarrow 0.$$ In addition,
$$P^N_{[e_0+\Delta e]}h=\frac {1}{\|e_0+\Delta e\|}f_{e_0+\Delta e}(h)(e_0+\Delta e) \quad P^{N_0}_{[e_0]}
h=\frac{1}{\|e_0\|}f_{e_0}(h)e_0 \quad \forall h\in E.$$ Obviously,
$e_0+\Delta e\rightarrow e_0$ as $\Delta e\rightarrow 0.$ Therefore,
one asserts $$f_{e_0+\Delta e}\rightarrow f_{e_0} \quad \rm{as}
\quad \Delta e\rightarrow 0.$$ Finally, one gets
$$(D\|\cdot\|)(e_0+\Delta e)=\|e_0+\Delta e\|f_{e_0+\Delta e}\rightarrow
 (D\|\cdot\|)(e_0)=\|e_0\|f_{e_0} \quad \rm{as} \ \Delta e \rightarrow 0.$$
 i.e., $(D\|\cdot\|)(e)$ is continuous at each $e_0 \in E\backslash \{0\}.$
\hfill $\Box$\\

Combining Theorems  2.1  and  3.1  bears the main theorem in the
paper:
\begin{th}\label {3.2} Suppose that $E$ is a Banach space. Then
the norm $\|\cdot\|$ of $E$ is of $c^1$ in $F\backslash \{0\}$ if and
only if $S(E)$ is a $c^1$-submanifold
 of $E$ with ${\rm codim} S(E)=1.$
\end{th}
\begin{co}\label{3.1} Suppose that $S$ is a $c^1$
submanifold of $E.$ Let $N=T_{e_0}S,  E=T_{e_0}S\oplus [e_0]$ for
$e_0\in S,$ and $e=P^{[e_0]}_Ne+\lambda e_0.$ then
$f_{e_0}(e)=\lambda \|e_0\| \in E^*$ fulfills $$(D\|\cdot\|)(e_0)h=f_{e_0}(h)
\quad \forall h\in E.$$
\end{co}
\section{Examples}\label{s4}
 \vspace{0.1cm}
\indent\quad The next two examples are interesting, which shows how to
determinate the Fr$\acute{e}$chet differential of the norm $\|\cdot\|$
by geometrical knowledge, although
they are simple.
\begin{ex}\label{1} Let $\|(x,y)\|=\sqrt{x^2+y^2}$ for any $(x,y)\in {\mathbb{R}}^2,$
 and $S$ be the unit circle with center 0. It is clear that the
 tangent line at a point $(x_0,y_0)\in S$ is the line perpendicular
 to the radial vector $(x_0,y_0),$ so that $N_0=T_{(x_0,y_0)}S=\{(x,y):
 xx_0+yy_0=0\}.$ Hence $${\mathbb{R}}^2=N_0\oplus [(x_0,y_0)],$$
and  $$P^{N_0}_{[(x_0,y_0)]}h=\lambda (x_0,y_0)=\lambda \|(x_0,y_0)\|e_0 $$
 for any $h\in {\mathbb{R}}^2,$
where $e_0=\frac {(x_0,y_0)}{\sqrt{x_0^2+y_0^2}}.$ As is well-known
from element geometry, reads the formula of the distance from $h$ to
 the tangent line of $S$ at $(x_0,y_0)$
$$f_{(x_0,y_0)}h=\lambda \|(x_0,y_0)\|=\frac{x_0h_1+y_0h_2}{\sqrt{x_0^2+y_0^2}}=x_0h_1+y_0h_2,$$
where $h=(h_1,h_2).$
By Theorem  3.1,
$$D(\sqrt{x^2+y^2})(x_0,y_0)h=x_0h_1+y_0h_2.$$
\end{ex}
\begin{ex}\label{2} Let $H$ be a Hilbert Space, $<,>$ denote its inner
product, and $\|h\|=\sqrt{<h,h>}.$ Let $S$ be the unit sphere in $H$
and $h_0\in S.$  Then  the subspace $N_0$  perpendicular to $h_0$ is just $T_{h_0}S$ and $N_0=T_{h_0}S=\{h\in H:
<H_0,h>=0.$ Since
 ${\rm codim}S=1$
 $$H= N_0\oplus [h_0].$$ Evidently,
$$P^{N_0}_{[h_0]}h=\lambda h_0=\lambda \|h_0\|e_0=<h, h_0>e_0 \quad \forall \ h\in H,$$
 where $e_0=\frac {h_0}{\|h_0\|}.$ By
Theorem  3.1,
$$(D\sqrt{<h,h>})(h_0)\Delta h=<h_0, \Delta h> \quad \forall \Delta
h\in H.$$
\end{ex}
 \vspace{0.2cm}
{\bf Acknowledgement}\quad  I would like to thank Professor
Yuwen Wang for telling me that he has considered the relationship between
smooth Banach space and smooth unit sphere for several years. This
work is supported by the National Science Foundation of China (Grant
No. 10671049 and 10771101).

1, Tseng Yuanrong Functional Research Center. Harbin Normal
University, Harbin 150080, P.R.China

2 Department of Mathematics, Nanjing University, Nanjing, 210093,
P.R.China\\
 E-mail address: \ jipuma@126.com

\end{document}